 
\iftrue
\else
\baselineskip=1.67\normalbaselineskip     
\smallskipamount=1.33\smallskipamount     
\medskipamount=1.33\medskipamount         
\bigskipamount=1.33\bigskipamount
\magnification=\magstep1
 
\hsize=5.9truein
\hoffset=0.5truein
\catcode`@=11
\def\plainoutput{\shipout\vbox{\makeheadline\pagebody\vskip10pt\makefootline}%
  \advancepageno
  \ifnum\outputpenalty>-\@MM \else\dosupereject\fi}
\catcode`@=12
\fi
 
\iftrue
\magnification=\magstep1
\hsize=6.5truein
\hoffset=0.0truein
\baselineskip 1.3\normalbaselineskip

\tolerance=10000
\def\sqr{$\vcenter{\hrule height .3mm
\hbox {\vrule width .3mm height 2mm \kern 2mm
\vrule width .3mm} \hrule height .3mm}$}
\else

\magnification=\magstep1
\baselineskip=1.33\normalbaselineskip
\tolerance=10000
\def\sqr{\vcenter{\hrule height .3mm
\hbox {\vrule width .3mm height 2mm \kern 2mm
\vrule width .3mm} \hrule height .3mm}}
\fi

\input amssym.def
\input amssym.tex
\def \chapter{\centerline}
\def \title{\medskip\centerline}

\def \Proof{\noindent {\bf Proof.\ \ }}

\def \E{{\Bbb E}}

\def \C{{\Bbb C}}

\def \ge{\geqslant}
\def \le{\leqslant}

\def \exp{\mathop{\rm exp}}

\def \e{\epsilon}

\def \d{\delta}

\def \a{\alpha}

\def \ra{\rightarrow}

\def \Z{{\Bbb Z}}
\def \C{{\Bbb C}}
\def \E{\mathop{{\Bbb E}}}

\centerline{\bf A UNIFORM SET WITH FEWER THAN EXPECTED ARITHMETIC}
\centerline{\bf PROGRESSIONS OF LENGTH 4.}
\medskip
\centerline{W. T. GOWERS}
\bigskip

\centerline{\it Centre for Mathematical Sciences, Wilberforce Road, Cambridge CB3 0WB, UK}
\bigskip

\centerline{Dedicated to Endre Szemer\'edi on the occasion of his 80th birthday.}
\bigskip

\noindent {\bf Abstract.} An example is constructed of a subset $A\subset\Z_N$ of density ${1\over 2}+o(1)$  such that all the non-trivial Fourier coefficients of the characteristic function of $A$ are very small, but if $x,d\in\Z_N$ are chosen uniformly at random, then the probability that $x,x+d,x+2d$ and $x+3d$ all belong to $A$ is at most ${1\over 16}-c$, where $c>0$ is an absolute constant.
\bigskip

\centerline {\bf 1. Introduction.}
\medskip

The purpose of this short note is to disprove a conjecture made in 
[3], which gave a proof of Szemer\'edi's theorem [5] that generalized
to longer progressions Roth's analytic approach [4] to the case of progressions of length 3. Let 
us very briefly explain the motivation for the conjecture, and the reason
that it seemed plausible: more details can be found in [3].

The main guiding principle behind the proof of [3], and indeed all 
other known proofs of Szemer\'edi's theorem, is that a set that is
sufficiently ``random-like'' contains about as many arithmetic progressions
of length $k$ as a typical random set of the same density. However, it
is not immediately obvious what should count as a ``random-like'' set.
If one is interested in progressions of length 3, then a good definition
turns out to come from Fourier analysis. First, one associates with a
subset $A\subset\{1,2,\dots,N\}$ its characteristic function. Following
[3], we shall write $A(x)$ rather than $\chi_A(x)$. Next, one regards
this function as a function defined on $\Z_N\equiv\Z/N\Z$ so that there
is a group structure appropriate for discrete Fourier analysis. Having
done so, one defines a discrete Fourier transform as follows. If $f:\Z_N\ra\C$
then 
$${\hat f}(r)=\E_{x\in\Z_N}f(x)\omega^{-rx}\ ,$$
where $\omega=\exp(2\pi i/N)$, and $\E$ is a standard abbreviation for $N^{-1}\sum$.

Note that ${\hat f}(0)=\E_xf(x)$, so if we return to our set $A$,
then ${\hat A}(0)=|A|/N$, or in other words the density of $A$. It turns out that an appropriate definition 
of quasirandomness for progressions of length three is that ${\hat A}(r)$
should be significantly smaller than this for all non-zero $r$. To be
precise, one can say that $A$ is a $c$-{\it uniform} set if 
$|{\hat A}(r)|\le c$ for every non-zero $r$.

A well-known and straightforward calculation shows that if $A$ has density $\a$
and is $c$-uniform for $c$ significantly smaller than $\a^2$, then 
$$\E_{x+z=2y}A(x)A(y)A(z)\approx\a^3.$$
That is, the ``3-AP density" of $A$ is roughly $\a^3$, which is roughly what it would
be (with high probability) for a random set of density $\a$. (Note that here we are 
counting mod-$N$ arithmetic progressions.) 

What makes it hard to prove Szemer\'edi's theorem analytically is that 
a set can be $c$-uniform for a very small $c$, even one that tends
to zero quite fast as $N$ tends to infinity, but yet not contain 
roughly the expected number of arithmetic progressions of length 4.
An example was given in [3], which we shall discuss in a moment.
However, a defect of the set defined in that example was that the 
number of progressions of length 4 was {\it greater} than one would
expect for a random set, so it did not demonstrate in a wholly 
satisfactory way that $c$-uniformity could not be used. At one
stage, I thought it would be routine to convert the example into
one for which the number of progressions was smaller, but eventually
Gil Kalai asked me how I proposed to do it, and I realized that it
was not routine after all. In fact, all the obvious ideas seemed
to fail, and the result was the conjecture in [3], which says
(in its mod-$N$ version) that if $A$ is a set of density $\d$ that 
is $c$-uniform, then the number of quadruples in $A^4$ of the form 
$(x,x+d,x+2d,x+3d)$ is at least $(\d^4-c')N^2$, where $c'$ tends
to zero as $c$ tends to zero. It is this conjecture that we shall
disprove here.

The example from [3] that yields too many progressions of 
length 4 is easy to define. Fix a small constant $c$, such as
$1/1000$, and let $A\subset\Z_N$ be the set of all $x\in\Z_N$
such that the least residue of $x^2$ mod $N$ lies in the interval 
$[-cN,cN]$. It is not too hard to verify that $A$ is 
$N^{-1/2+\e}$-uniform (which, since $\sum_r|{\hat A}(r)|^2\le 1$ 
is almost as uniform as a set can possibly be) and has cardinality
approximately $2cN$: proofs of these facts will be sketched later 
in this paper. Therefore, by the argument above, $A$ contains about 
$8c^3N^2$ triples $(x,x+d,x+2d)$. In fact, one can say a little
more: triples of the form $(x^2,(x+d)^2,(x+2d)^2)$ are ``approximately
uniformly distributed'' in the sense that, given any three reasonably
long intervals $I$, $J$ and $K$, they contain about $N^{-1}|I||J||K|$
of them. Now we use this uniform distribution property, together 
with the identity 
$$(x+3d)^2=x^2-3(x+d)^2+3(x+2d)^2$$
to conclude that, if $x$ and $d$ are chosen randomly, then the 
probability that $(x+3d)^2$ belongs to $[-cN,cN]$ given that $x^2$, 
$(x+d)^2$ and $(x+2d)^2$ all belong to $[-cN,cN]$ is significantly 
bigger than $2c$ (indeed, it will be approximately equal to an
absolute constant that one could work out after a tedious calculation), 
and therefore that the number of quadruples $(x,x+d,x+2d,x+3d)$ in
$A^4$ is significantly greater than $16c^4N^2$.

A natural way to try to generalize this example is by replacing the
interval $[-cN,cN]$ with a more general set $B\subset\Z_N$. For the above 
argument to work, we would want $B$ to be a union of longish intervals, so
that the uniform distribution property was still valid. Then we 
could define $A$ to be $\{x\in\Z_N:x^2\in B\}$. If $B$ has size
$\beta N$, then the uniform distribution property implies that
the number of triples in $A$ of the form $(x,x+d,x+2d)$ is approximately
$\beta^3N^2$, and that the corresponding triples $(x^2,(x+d)^2,(x+2d)^2)$ 
are ``approximately uniformly distributed'' in $B^3$. We now want to
estimate the conditional probability that $(x+3d)^2\in B$, given that
$x^2$, $(x+d)^2$ and $(x+2d)^2$ all lie in $B$. 

If the distribution were genuinely uniform, then, by the identity we 
used earlier, this conditional probability would be $|B|^{-3}$ times
the number of quadruples $(x,y,z,w)\in B^4$ such that $w=x-3y+3z$. The
notion of uniform distribution turns out to be strong enough for this
to be a good estimate even in the situation we actually face, so let
us think about how many such quadruples there are, by expressing the
number in terms of Fourier coefficients. We obtain
$$\eqalign{\E_{x-3y+3z-w=0}B(x)B(y)B(z)B(w)&=
\sum_r\sum_{x,y,z,w}B(x)B(y)B(z)B(w)\omega^{r(x-3y+3z-w)}\cr
&=\sum_r|{\hat B}(r)|^2|{\hat B}(3r)|^2\cr
&=\beta^4+\sum_{r\ne 0}|{\hat B}(r)|^2|{\hat B}(3r)|^2\,\cr}$$
where $\beta=|B|/N$ is the density of $B$. 
Therefore, the conditional probability
in question is at least $\beta$, since the contribution from the non-zero
$r$ is positive. Therefore, the number of progressions of length 4 in $A$ 
is bounded below (at least approximately) by $\beta^4N^2$. 

It is not necessary to use Fourier analysis to obtain this
conclusion: one can also observe that the left-hand side is the square of
the $L_2$ norm of the convolution of $B$ with a dilate of $B$, and apply
Cauchy-Schwarz. In any case, the conclusion shows that no $A$ constructed in this 
way can give an example of a set with too few progressions of length 4. 

Note that the above argument depends on the positivity of the Fourier
expression, which in turn depends on the fact that the set $1,-3,3,-1$ 
can be partitioned into pairs of the form $\{a,-a\}$. The non-Fourier
argument also relies on this fact, in a slightly different way. Thus, the fact
that 4 is an even number is highly relevant. Indeed, this obstruction
does not exist for odd $k$: it is possible to choose a union of intervals
$B$ of density $\beta$ such that the number of quintuples $(v,w,x,y,z)\in B^5$ 
with $v-4w+6x-4y+z=0$ is significantly less than $\beta^5N^5$. Then one
can define $A$ to be $\{x\in\Z_N:x^3\in B\}$. 

This set can be shown to have small Fourier coefficients. It turns out that the
uniform distribution property applies to quadruples of the form
$(x^3,(x+d)^3,(x+2d)^3,(x+3d)^3)$, from which it follows that $A$
contains fewer than $c\beta^5N^2$ quintuples of the form
$(x,x+d,x+2d,x+3d,x+4d)$, where $c<1$ is an absolute constant. The details 
of this argument were worked out
in a conversation with Ben Green, but here we leave them as an exercise
for the interested reader and return to progressions of length 4.

We have now given the main motivation for the false conjecture in [3]:
there seems to be a clear difference between 4 and 5, arising from their
differing parities; for progressions of length 4, there is a certain
positivity phenomenon that stops examples of a certain kind from having
too few progressions of length 4. However, with a little more ingenuity 
one can after all produce an example that works for progressions of 
length 4, as we shall now show. 
\bigskip

\centerline {\bf 2. Construction of a uniform set with few progressions 
of length 4.}
\medskip

\proclaim Lemma 1. There is a function $f:\Z\ra\{-1,0,1\}$ such that
$f(n)=0$ unless $1\le n\le 300$ and such that 
$$\sum_{x,d\in\Z}f(x)f(x+d)f(x+2d)f(x+3d)<0\ .$$

\Proof  First let us define a function $g:\Z^3\ra\{-1,0,1\}$ 
with the required property (where now $x$ and $d$ are elements of
$\Z^3$) that is supported in the set $\{1,2,3,4\}^3$. Given
$(a,b,c)\in\{1,2,3,4\}^3$ we let $g(a,b,c)=1$ unless $(a,b,c)$ is one
of the following sixteen triples, in which case we set $g(a,b,c)=-1$:
113,121,132,144,212,224, 233,241,314,322,331,343,411,423,434 and
442. These triples are chosen to have the following property: every
line that goes through four points of the grid $\{1,2,3,4\}^3$ and is
not one of the four main diagonals contains precisely one of these
points. (The search for such a system of points was not very
difficult, since in each horizontal plane one had to take exactly one
point in each row and column, and this had to be done in a disjoint
way for the four planes. This was already enough of a restriction to
make aesthetically driven trial and error a feasible method.)

Each pair $(x,d)$ that contributes to the sum above, with $g$
instead of $f$, is either
degenerate, in the sense that $d=0$, or it corresponds to a
geometrical line. Note also that each non-degenerate line is
counted twice, since the line for $x$ and $d$ is the same as
the line for $x+3d$ and $-d$. The contribution to the sum from
each degenerate pairs is $1$, and there are $4^3=64$ of them.
The contribution from each line that is parallel to one of the
three coordinate axes is $-2$ and there are $3\times 16=48$ of
these. The contribution from each non-main diagonal is $-2$
and there are $3\times 8=24$ of those. Finally, each of the
four main diagonals contributes $2$. Therefore, the sum is
$64-96-48+8=-72$.

Now we transfer this example to $\Z$ by a suitable projection.
Define a map $\phi:\{1,2,3,4\}^3\ra\Z$ by 
$\phi(a,b,c)=a+8b+64c$. This map is easily checked to be a 
Freiman homomorphism: that is, $\phi(x)-\phi(y)=\phi(z)-\phi(w)$
if and only if $x-y=z-w$. Indeed, if
$$a_1+8b_1+64c_1-a_2-8b_2-64c_2=a_3+8b_3+64c_3-a_4-8b_4-64c_4\ ,$$
then 
$$(a_1-a_2-a_3+a_4)+8(b_1-b_2-b_3+b_4)+64(c_1-c_2-c_3+c_4)=0\ .$$
But $a_1-a_2-a_3+a_4$, $b_1-b_2-b_3+b_4$ and $c_1-c_2-c_3+c_4$
all lie between $-6$ and $6$, from which it follows easily that
they are all zero.

A quadruple of the form $(x,x+d,x+2d,x+3d)$ is the same as a
quadruple $(x,y,z,w)$ such that $y-x=z-y$ and $z-y=w-z$. Therefore,
$(\phi x, \phi y, \phi z, \phi w)$ is such a quadruple if and
only if $(x,y,z,w)$ is. It follows that if we define
$f(\phi(a,b,c))$ to be $g(a,b,c)$ and $f(x)$ to be 0 if 
$x$ is not in the image of $\phi$, then 
$$\sum_{x,d}f(x)f(x+d)f(x+2d)f(x+3d)=-72\ .$$
This proves the lemma. \hfill $\square$ \bigskip

\noindent {\bf Remark.} A more efficient version of the above lemma was 
obtained by Wolf [6], who found a function $f$ supported on $\{1,2,\dots,18\}$ 
that takes only $\pm 1$ values in that range, for which the sum is -36. 
\bigskip

\proclaim Corollary 2. There exists an absolute constant $c>0$ such 
that for every sufficiently large $N$ there is a function
$F:\Z_N\ra\{-1,0,1\}$ with the following properties. First,
$F$ is a $\pm 1$-combination of characteristic functions of
disjoint intervals of size at least $N/1500$, and, secondly,
$$\E_{x,d\in\Z_N}F(x)F(x+d)F(x+2d)F(x+3d)\le -c\ .$$

\Proof  Let $t$ be a positive integer between $N/1500$ and
$N/1200$, and for $1\le k\le 300$ let $I_k$ be the interval
$\{(2k-1)t+1,(2k-1)t+2,\dots,2kt\}$. Then an argument similar
to the proof that the function $\phi$ in Lemma 1 was a Freiman 
homomorphism shows that if $x$ and $d$ are elements of $\Z_N$
and $x\in I_k$, $x+d\in I_l$, $x+2d\in I_m$ and $x+3d\in I_n$,
then $l-k=m-l=n-m$. In other words, an arithmetic progression
can lie in the union of the intervals $I_k$ only if the 
corresponding intervals themselves lie in an arithmetic 
progression (which may be degenerate).

Suppose that $(k,l,m,n)$ is an arithmetic progression. Then
the number of arithmetic progressions 
$(x,y,z,w)\in I_k\times I_l\times I_m\times I_n$ is the same
as the number of such progressions in $\{1,2,\dots,t\}^4$,
since the map 
$$(x,y,z,w)\mapsto (x-(2k-1)t,y-(2l-1)t,z-(2m-1)t,w-(2n-1)t)$$
is a bijection between the two sets of progressions in question.
Notice also that every mod-$N$ progression that lives in the 
union of the $I_k$ must in fact correspond to a genuine 
arithmetic progression in the set $\{1,2,\dots,N\}$, since
the union does not contain any numbers between $N/2$ and $N$.

Let $p$ be the number of progressions in $\{1,2,\dots,t\}^4$, where
by ``progression'' we mean ``quadruple of the
form $(x,x+d,x+2d,x+3d)$'' and allow $d$ to be any element
of $\Z_N$. Then $p\ge t^2/3$, since if we choose
two points at random in $\{1,2,\dots,t\}$ there is a probability
at least ${1\over 3}$ that they differ by a multiple of 3.

Now let $F$ be defined as follows. If $x\in I_k$ then $F(x)=f(k)$.
For all other $x$, $F(x)=0$. The remarks we have just made 
show that
$$\E_{x,d\in\Z_N}F(x)F(x+d)F(x+2d)F(x+3d)
=pN^{-2}\sum_{x,d\in\Z}f(x)f(x+d)f(x+2d)f(x+3d)\ .$$
But $p\ge t^2/3$, $t\ge N/1500$ and 
$\sum_{x,d\in\Z}f(x)f(x+d)f(x+2d)f(x+3d)=-72$. Therefore,
the lemma is proved, with $c=72/3(1500)^2$, which is 
greater than $1/100000$. \hfill $\square$ \bigskip

Corollary 2 may seem rather pointless, since the function $F$ is
composed of characteristic functions of intervals and therefore
has large Fourier coefficients. However, the main observation that
underlies the construction is that such an example can be converted
into a uniform set in a rather simple way, while keeping its
property of having a negative sum over progressions of length 4.

The trick is to multiply $F$ pointwise by a small linear combination
of ``quadratic phase functions,'' that is, functions of the form
$\omega^{q(x)}$ for a quadratic polynomial $q$. These have two
properties that make them useful. First, quadratic functions 
satisfy the identity 
$$q(x)-3q(x+d)+3q(x+2d)-q(x+3d)=0\ ,$$
as we have already noted. Secondly, functions of the form 
$\omega^{q(x)}$ are highly uniform -- indeed, so uniform that
they remain highly uniform even when multiplied by characteristic
functions of intervals. Let us prove the second of these facts,
which will allow us to show that the function we create is uniform.
The proof is standard, but is included for convenience. Once we have
established the second fact, we can show why the first fact helps 
us to obtain a good estimate for sums of products over arithmetic 
progressions.

\proclaim Lemma 3. Let $I$ be an interval of length $t$ in $\Z_N$
and let $q$ be a quadratic polynomial. Then no Fourier coefficient
of the product $I(x)\omega^{q(x)}$ is greater than $2N^{-1/2}\log N$.

\Proof  First we show that all Fourier coefficients of the function 
$Q(x)=\omega^{q(x)}$ have size $N^{-1/2}$. Indeed, let 
$q(x)=ax^2+bx+c$, with $a\ne 0$. Then
$$\eqalign{\Bigl|\E_xQ(x)\omega^{-rx}\Bigr|^2
&=\E_{x,y}\omega^{-(q(x)-q(y)+r(x-y))}\cr
&=\E_{x,y}\omega^{-(ax+ay+b+r)(x-y)}\cr
&=\E_z\E_x\omega^{-(2ax-az+b+r)z}\ .\cr}$$
For each fixed $z$, the average over $x$ is zero, except when 
$z=0$, in which case it is $1$. Therefore, the average equals
$N^{-1}$, which proves the claim.

We now use the convolution identity in the less standard direction,
which tells us that the $r$th Fourier coefficient of the pointwise product
$QI$ is $\sum_{s+t=r}\hat Q(s)\hat I(r)$, which has modulus at most
$N^{-1/2}\sum_r|{\hat I}(r)|$.

By the formula for summing a geometric progression, we find that
$|\hat I(r)|$ is at most $2/N|1-\omega^r|$. One can also check that 
$|1-\omega^r|\geq 4r/N$ when $-N/2\leq r\leq N/2$. It follows that 
$\sum_r|{\hat I}(r)|$ is at most $1+\sum_{r\le N/2}r^{-1}$,
which is at most $2\log N$. The result follows.~\hfill~$\square$
\bigskip

Now let us define a function $G$ by
$$G(x)=F(x)(\omega^{x^2}+\omega^{-x^2}+\omega^{3x^2}+\omega^{-3x^2})\ .$$
Then, since $F$ is a $\pm 1$-combination of characteristic functions
of $64$ intervals, Lemma 3 implies that no Fourier coefficient of
$G$ exceeds $4\times 64\times 2N^{-1/2}\log N=512N^{-1/2}\log N$. Notice 
also that $G$ is real valued and that its values all belong to the 
interval $[-4,4]$.

\proclaim Lemma 4. Let $G$ be as just defined. Then
$$\Bigl|\E_{x,d}G(x)G(x+d)G(x+2d)G(x+3d)
-2\E_{x,d}F(x)F(x+d)F(x+2d)F(x+3d)\Bigr|$$
is at most $2^{18}N^{-1/2}\log N$.

\Proof The product $G(x)G(x+d)G(x+2d)G(x+3d)$
splits into a sum of $256$ terms of the form
$\omega^{px^2+q(x+d)^2+r(x+2d)^2+s(x+3d)^2}$, where $p$, $q$, $r$
and $s$ belong to the set $\{-3,-1,1,3\}$. Let us fix a choice
of $p$, $q$, $r$ and $s$, set $\theta(x,d)=px^2+q(x+d)^2+r(x+2d)^2+s(x+3d)^2$,
and estimate the quantity 
$$\E_{x,d}F(x)F(x+d)F(x+2d)F(x+3d)\omega^{\theta(x,d)}\ .$$
We can write $\theta(x,d)=ux^2+vxd+wd^2$. Let us suppose that $u\ne 0$
and consider the average for a fixed value of $d$, writing 
$F_d(x)$ for $F(x)F(x+d)F(x+2d)F(x+3d)$. The number of intervals
$I_k$ on which $F$ is non-zero is 64, so as $x$ increases, the number
of times at least one of $x$, $x+d$, $x+2d$ and $x+3d$ changes from
not belonging to a certain $I_k$ to belonging to it is at most
$4\times 64=256$. It follows that $F_d$ is a $\pm 1$-sum of
characteristic functions of at most $256$ disjoint intervals.
Lemma 3 (in the case of the Fourier coefficient at zero) then
implies that $\E_xF_d(x)\omega^{\theta(x,d)}$ has modulus at
most $512N^{-1/2}\log N$. If we now average over $d$ we find that
$$\E_{x,d}F(x)F(x+d)F(x+2d)F(x+3d)\omega^{\theta(x,d)}\le 512N^{-1/2}\log N\ .$$

Now let us consider what happens if $w\ne 0$. The argument is
very similar. This time we shall fix $x$ and define $F_x(d)$ to
be $F(x)F(x+d)F(x+2d)F(x+3d)$. (Note that we are using the
same notation for a different definition in this case.) As
$d$ increases, the number of times $x+rd$ can change from not
belonging to $I_k$ to belonging to $I_k$ is at most $r$. Therefore,
$F_x$ is a $\pm 1$-sum of characteristic functions of at most 
$(1+2+3)\times 64=384$ disjoint intervals. So in this case
we obtain the estimate
$$\E_{x,d}F(x)F(x+d)F(x+2d)F(x+3d)\omega^{\theta(x,d)}\le 768N^{-1/2}\log N\ .$$

If $u=0$ and $w=0$ then $p+q+r+s=0$ and $q+4r+9s=0$. Because
$p$, $q$, $r$ and $s$ all lie in the set $\{-3,-1,1,3\}$, it
is not possible for $q+4r+9s$ to equal 0 unless $s=\pm 1$. 
Since we can multiply any solution by $-1$ and obtain another
solution, let us suppose that $s=-1$. Then $q+4r=9$, which, it
is simple to check, can happen only if $q=-3$ and $r=3$. But
then $p=1$. It follows that the only two solutions are
$(1,-3,3,-1)$ and $(-1,3,-3,1)$. 

But in this case, $\theta(x,d)=0$, by the identity used earlier,
so the average becomes $\E_{x,d}F(x)F(x+d)F(x+2d)F(x+3d)$. It
follows that the difference between the sum with $G$ and
the sum with $F$ is at most $256\times 1024N^{-1/2}\log N$,
which is what the lemma states.~\hfill~$\square$ \bigskip

We are almost done. $G$ is a function taking values in 
$[-4,4]$ with a ``negative number of progressions of
length 4.'' All that remains is to convert it into a set.
This we do by first setting $P(x)$ to be $(G(x)+4)/8$
and then choosing a set $A$ by letting $x$ belong to
$A$ with probability $P(x)$. Here are the details.

\proclaim Corollary 5. There exists an absolute constant
$c>0$ such that for every sufficiently large $N$ there is
a function $P:\Z_N\ra[0,1]$ such that
$\Bigl|\E_xP(x)-{1\over 2}\Bigr|\le 64N^{-1/2}\log N$, such
that $|{\hat P}(r)|\le 64N^{-1/2}\log N$ for every $r\ne 0$
and such that 
$$\E_{x,d}P(x)P(x+d)P(x+2d)P(x+3d)\le {1\over 16}-c\ .$$

\Proof  As we have already said, we let $P(x)$ equal
$(G(x)+4)/8$. Then $|\E_xP(x)-{1\over 2}|=|\E_xG(x)/8|$, which
is at most $64N^{-1/2}\log N$, since it is ${\hat G}(0)/8$
and all Fourier coefficients of $G$ are at most 
$512N^{-1/2}\log N$. Similarly, when $r\ne 0$,
${\hat P}(r)={\hat G}(r)/8$, since adding a constant
to a function does not alter its Fourier coefficients
for $r\ne 0$.

To prove the last property, write the product
$P(x)P(x+d)P(x+2d)P(x+3d)$ as
$$2^{-12}(4+G(x))(4+G(x+d))(4+G(x+2d))(4+G(x+3d))\ .$$
This product splits into 16 parts, each of which we
shall average separately.

If we choose $4$ from every bracket, then we obtain 
$2^{8-12}={1\over 16}$ for every $x$ and $d$. 
If we choose $G$ from every bracket then
we are estimating 
$$2^{-12}\E_{x,d}G(x)G(x+d)G(x+2d)G(x+3d)\ ,$$
which, by Corollary 2 (with the quantitative estimate
at the end) and Lemma 4 is at most $-2^{-12}.10^{-5}+2^6N^{-1/2}\log N$.

If we choose $G$ from precisely one bracket, then we
obtain $2^{-6}\E_xG(x)$, after a suitable change of
variables. We have already remarked that 
$\left|\E_xG(x)\right|\le 512N^{-1/2}\log N$, so
this is at most $8N^{-1/2}\log N$.

If we choose $G$ from precisely two brackets, then
we obtain $2^{-8}\left(\E_xG(x)\right)^2$, again 
after a suitable change of variables. This has modulus
at most $1024N^{-1}(\log N)^2$.

If we choose $G$ from three brackets, then we must
do four standard calculations similar to the one that proves that
a uniform set contains many arithmetic progressions
of length 3. We illustrate this by bounding the quantity
$$2^{-10}\E_{x,d}G(x)G(x+d)G(x+3d)\ .$$
First, a routine Fourier calculation shows that it is equal to
$$\sum_r\hat G(2r)\hat G(-3r)\hat G(r)\ .$$
But $|{\hat G}(r)|\le 512N^{-1/2}\log N$ for every $r$. Using
this fact and the Cauchy-Schwarz inequality, the last expression
can be bounded above in modulus by
$$512N^{-1/2}\log N\Bigl(\sum_r|{\hat G}(2r)|^2\Bigr)^{1/2}
\Bigl(\sum_r|{\hat G}(-3r)|^2\Bigr)^{1/2}\ .$$
Since $N$ is prime, the product of the last two
brackets is $\|\hat G\|_2^2=\|G\|_2^2\leq 16$.
Therefore, the
whole sum is at most $2^{13}N^{-1/2}\log N$. After multiplying
by $2^{-10}$, as we need to, we have shown that the contribution 
from this term is at most $8N^{-1/2}\log N$. The estimates for
the other terms with $G$ chosen three times are proved in a
very similar way.

There are sixteen terms, and the worst error term (when $N$ is
sufficiently large) has modulus at most $8N^{-1/2}\log N$.
Therefore, $\E_{x,d}P(x)P(x+d)P(x+2d)P(x+3d)$ is at 
most $(2^{-4}-2^{-12}.10^{-5})+(2^6+2^7)N^{-1/2}\log N$. For
sufficiently large $N$, this is less than ${1\over 16}-2^{-29}$, which proves the corollary.
\hfill $\square$ \bigskip

\proclaim Theorem 6. There exist absolute constants
$c>0$ and $C$ such that for all sufficiently large $N$ there is
a $CN^{-1/2}\log N$-uniform subset $A$ of $\Z_N$ of density 
${1\over 2}+o(1)$, for which 
$$\E_{x,d}A(x)A(x+d)A(x+2d)A(x+3d)\le 2^{-4}-c\ .$$

\Proof  Choose $A$ randomly as follows. For every $x\in\Z_N$, 
let $x$ belong to $A$ with probability $P(x)$, with all 
choices made independently. Then $A(x)-P(x)$ is a random
variable of mean zero that is equal to either $1-P(x)$
or $-P(x)$. In particular, it always has modulus at most
1. 

Now let us consider the difference between the Fourier
coefficients of $A$ and those of $P$. We have
$${\hat A}(r)-{\hat P}(r)=\E_x(A(x)-P(x))\omega^{rx}\ .$$
This is an average of $N$ independent random variables, each of mean zero
and bounded above in modulus by 1. It follows immediately
from Azuma's inequality that the probability
that $|{\hat A}(r)-{\hat P}(r)|$ is greater than $tN^{-1/2}$
is at most $2e^{-t^2/8}$. If we take $t$ to equal $\log N$,
then this probability is, for $N$ sufficiently large, 
significantly smaller than $1/N^2$. It follows that with
probability at least $1-1/N$ we have 
$|{\hat A}(r)-{\hat P}(r)|\le N^{-1/2}\log N$ for every
$r$. This implies both the statement about the cardinality
of $A$ and the claim about its uniformity.

For each non-zero $d$ the expectation of the
quantity $\E_xA(x)A(x+d)A(x+2d)A(x+3d)$ is equal to
$\E_xP(x)P(x+d)P(x+2d)P(x+3d)$, and when $d=0$
it is at most $1$. Therefore, the expectation of
$\E_xA(x)A(x+d)A(x+2d)A(x+3d)$ is at most
${1\over 16}-c+N^{-1}$, where $c$ is the absolute constant
obtained in Lemma 5. Therefore, the probability that
$\E_xA(x)A(x+d)A(x+2d)A(x+3d)$ is at most
${1\over 16}-{c\over 2}$ is at least $8c$. Since this is
bigger than $N^{-1}$ when $N$ is sufficiently large,
we can find a choice of $A$ that has all the
properties claimed. \hfill $\square$ \bigskip

The main result of this paper answers the conjecture
in [3] but it does not answer every question
one might wish to ask. In particular, it is noticeable
that the number of progressions of length 4 in $A$ was
not {\it much} less than it would be for a random set.
The following question, which was asked by  
Ruzsa and appears in [2] as Problem 3.2, is still open.

\proclaim Problem. Let $A$ be a uniform subset of $\Z_N$
of density $\a$. Must $A$ contain at least $\a^{1000}N^2$
arithmetic progressions of length 4?
\bigskip

\noindent Of course ``1000" is just an informal way of 
referring to {\it some} absolute constant. This question is open 
only for progressions of length 4: for longer 
progressions, Ruzsa has constructed a counterexample (see Theorem 2.4
of [1]).
\bigskip

\centerline {\bf References.}
\medskip

\noindent [1] V. Bergelson, B. Host. and B. Kra. with an appendix by I. Z. Ruzsa, Multiple recurrence and nilsequences, {\it Invent. Math.}, {\bf 160} (2005), 261-303.
\medskip

\noindent [2] E. S. Croot and V. F. Lev, Open problems in additive combinatorics, in: {\it Additive Combinatorics}, CRM Proc. Lecture Notes 43, Amer. Math. Soc. (Providence, RI 2007), pp. 207-233.
\medskip

\noindent [3] W. T. Gowers, A new proof of Szemer\'edi's theorem, {\it Geom. Funct. Anal.} {\bf 11} (2001), 465-588. 
\medskip

\noindent [4] K. F. Roth, On certain sets of integers, {\it J. London Math. Soc.} {\bf 28} (1953), 245-252.
\medskip

\noindent [5] E. Szemer\'edi, On sets of integers containing no $k$ elements in arithmetic progression, {\it Acta Arith.} {\bf 27} (1975), 299-345.
\medskip

\noindent [6] J. Wolf, The minimum number of monochromatic 4-term progressions in $\Bbb Z_p$, {\it J. Comb.} {\bf 1} (2010), 53-68.

\bye